\documentclass[%
preprint,
 amsmath,amssymb,
 aps,
]{revtex4-1}

\usepackage{amsmath,mathtools,amsthm, tikz,amssymb,enumerate,framed,listings,multirow,pdftricks,subcaption,times,xcolor,soul,booktabs,graphicx,epstopdf}

\usepackage{ifthen}

\colorlet{shadecolor}{gray!0}

\def\Y_#1{\vect{y}_{\!#1}}

\def\E{\mathbb{E}}

\def\tr{\operatorname{trace}}

\def\spanset{\operatorname{span}}

\newcommand{\vect}[1]{\boldsymbol{#1}}

\newtheorem{theorem}{Theorem}[section]
\newtheorem{thm}{Theorem}[section]

\theoremstyle{definition}

\newtheorem{algorithm}[theorem]{Algorithm}

\theoremstyle{remark}
\newtheorem{remark}[theorem]{Remark}

\def\E{\mathbb{E}}

\begin{document}

\lstset{language=Matlab}          
\lstset{  basicstyle=\footnotesize,        
  breakatwhitespace=false,         
  breaklines=false,                 
numbers=left, numberstyle=\tiny, stepnumber=1, numbersep=5pt, keywordstyle=\bfseries,basicstyle=\ttfamily,  frame=lines}

\title{
SPECTRWM: Spectral Random Walk Method  for the \\
Numerical Solution of Stochastic Partial Differential Equations
}
\author{Nawaf Bou-Rabee}
 \email{nawaf.bourabee@rutgers.edu}
\affiliation{
Department of Mathematical Sciences  \\
Rutgers University \\
311 N 5th Street \\
Camden, NJ 08102, USA
}

\begin{abstract}
The numerical solution of stochastic partial differential equations (SPDE) presents challenges not encountered in the simulation of PDEs or SDEs.  Indeed, the roughness of the noise in conjunction with nonlinearities in the drift typically make these equations particularly stiff.  In practice, this means that it is tricky to construct, operate, and validate numerical methods for SPDEs.  This is especially true if one is interested in path-dependent expected values, long-time simulations, or in the simulation of SPDEs whose solutions have constraints on their domains.  To address these numerical issues, this paper introduces a Markov jump process approximation for SPDEs, which we refer to as the spectral random walk method (SPECTRWM).   The accuracy and ergodicity of SPECTRWM are verified in the context of a heat and overdamped Langevin SPDE, respectively.  We also apply the method to Burgers and KPZ SPDEs.  
\end{abstract}

%
%
\maketitle

\section{Introduction}

Stochastic partial differential equations (SPDEs) describe the evolution of continuum systems with random fluctuations \cite{Wa1986, Ha2009, DaZa2014}.  They are used as models for turbulence \cite{WeKhMaSi1997, WeVa2000, EKhMaSi2003, MiRo2004, DaDe2003, HaMa2006}, phase-field dynamics \cite{Fu1995, KoOtReVa2007,We2010},  surface growth \cite{KPZ1986, Qu2011, Co2012, Ha2013}, neuronal activity \cite{Wa1986, Tu2005, LaLo2010}, population dynamics/genetics \cite{Fe1951,Da1972,Fl1975} and interest rate fluctuations \cite{Co2005, CaTe2007}.  They also directly arise as space-time diffusion approximations of a large class of discrete models \cite{KPZ1986, Qu2011, Co2012, DeTs2015}.  In addition to their use as mathematical models, their dynamic and ergodic properties are also leveraged to accelerate convergence of MCMC methods for sampling conditioned diffusions \cite{StVoWi2004, ReVa2005, HaStVoWi2005, HaStVo2007, ApHaStVo2007, BeRoStVo2008}.   The importance of these equations lies in the basic fact that space-time Gaussian white noise is perhaps the simplest model of space-time random fluctuations.  In this sense, the study of SPDEs is a quest to answer fundamental questions about continuum descriptions of high-dimensional discrete systems.

With few exceptions \cite{HoLuRoZh2006,Lu2006}, current numerical methods for SPDEs are rooted in ideas from numerical PDE and SDE theory: replace the SPDE with a system of approximating SDEs by means of a finite difference or variational method, and then discretize these approximating SDEs in time using, e.g., a $\theta$-method like Euler-Maruyama or Crank Nicolson; see, e.g., \cite{GyNu1995, GyNu1997, DaGa2001, Wa2005, AlGy2006, JeKl2009, BlJe2013, Kr2014}.   
However, there are particular problems associated with approximating SPDEs that motivate treating them differently.  Here is a partial list.
\begin{itemize}
\item {\em Stiffness.}  
Due to the lack of regularity of the SPDE solution, the optimal time discretization error associated with spatio-temporal approximations is of strong order $1/4$ and weak order $1/2$ \cite{DaGa2001, Wa2005}.  (To obtain these rates, one must assume that the spatial step size relates to the time step size via a CFL-type condition.) These rates contrast with the faster order of $2$ typically attainable for spatio-temporal approximations of the underlying PDE, and of strong order $1/2$ and weak order of $1$ typically attainable for SDEs using the Euler-Maruyama method \cite{KlPl1995,Ta1995,MiTr2004}.  
\item {\em Spurious Drift Terms.}
Hairer et al.~showed that seemingly reasonable spatial discretizations of Burgers-like SPDEs do not converge to the right solution -- even in a weak sense \cite{HaVo2011, HaMa2012, HaMaWe2014}.   Instead these discretizations converge to a (non-ergodic) SPDE with a spurious drift term that is a spatial analog of the It\^o-Stratonovich correction.  This numerical artifact is due to the spatial roughness of the SPDE solution.   We refer the reader to \cite{HaVo2011} for numerical results and an explanation of this numerical artifact, and to \cite{HaMa2012, HaMaWe2014} for detailed proofs.
\item {\em Long-Time Simulation.}
Another problem is related to long-time simulation of ergodic SPDEs \cite{DaZa1996, cerrai2001second, Vo2011}.  The aim of this type of simulation is to sample from the stationary distribution of the SPDE, and also to compute long-time dynamics.  As far as we can tell, methods for ergodic SPDEs are currently limited to schemes that efficiently sample from the stationary distribution of Langevin SPDEs without too much concern about accurately representing their dynamics \cite{BeRoStVo2008}.  The idea for this scheme comes from MCMC and numerical SDE theory \cite{RoTw1996A, RoTw1996B, BoVa2010,BoVa2012,BoHa2013,BoDoVa2014, Bo2014, Fa2014, FaHoSt2015}, and basically involves combining a $\theta$-method for the approximating SDE with a Metropolis accept-reject step, which corrects the bias introduced by time discretization error.   However, the results in Ref.~\cite{BeRoStVo2008} demonstrate that unless one chooses $\theta=1/2$ (a Crank-Nicholson discretization), the acceptance rate of the Metropolized $\theta$-scheme deteriorates in the space-time diffusion limit.
\item {\em Reflecting Solutions.}
Reflected SPDEs are an infinite-dimensional analog of SDEs with reflection \cite{KaSh1991,NuPa1992,KuDu2001}.  They arise in the study of continuum models with constraints on the domain of their solutions, e.g., the Allen-Cahn SPDE with reflection at $\pm 1$ or positivity constraints in population or interest rate models \cite{NuPa1992,FuOl2001,DeZa2007,DeGo2011}.   In this context spatio-temporal methods may produce approximations that are outside the domain of definition of the SPDE.  
\end{itemize}

In light of these practical problems, it is quite natural to approximate SPDEs using an algorithm that is more precisely tailored to the structure of their solutions.   In this note we propose a simple way to do this.  The idea is to approximate SPDEs  by a Markov jump process that we refer to as the spectral random walk method (SPECTRWM).  The departure point for constructing this approximation is a stable and accurate system of approximating SDEs, which may be obtained by, e.g., finite difference or spectral Galerkin methods \cite{GyNu1995, GyNu1997, Wa2005}.  To be clear, we do not propose to solve the above problem regarding spurious drift terms, though we do confirm that incorrect discretizations lead to non-ergodic approximations to the Burgers and KPZ SPDEs.   Given these approximating SDEs, we then proceed as follows: instead of discretizing these approximating SDEs in time -- as is normally done -- we discretize their infinitesimal generator in space \cite{BoVa2016}.  As we detail in \S\ref{sec:spectrwm} below, an important ingredient in this construction is a basis given by the leading $n$ eigenfunctions associated to the linear part of the drift of the SPDE.

A realization of this Markov jump process approximation may be produced by iteratively computing its jumps and holding times:  the jumps are taken in the direction of these eigenfunctions, while the holding time in any given state is an exponentially distributed random variable whose mean is a deterministic function of the SPDE coefficients evaluated at this state.  In addition to being simple, the method allows several benefits: 
\begin{itemize}
\item the jump size of the approximation is a parameter of the method;
\item the time step size automatically adapts according to the stiffness of the SPDE coefficients; 
\item path-dependent expected values over slices of the SPDE solution in time can be approximated without incurring any time discretization error; 
\item every jump induces a global move in state space; 
\item they are multi-scale, in the sense that its jump size can be adapted to the different spatial scales of the SPDE problem; and,
\item it can handle boundary conditions in reflected SPDEs in a natural way.
\end{itemize}
The SPECTRWM method and its properties are described in \S\ref{sec:spectrwm}.  Afterwards we test SPECTRWM on heat, overdamped Langevin, Burgers, and KPZ SPDEs, all with periodic boundary conditions.

Let us finish this introduction by remarking that SPECTRWM is a generalization of the Markov Chain Approximation Method (MCAM) to SPDEs.  By now, MCAM is a well-established technique for the numerical solution of SDEs \cite{KuDu2001, BoVa2016}.  The method was invented by Harold Kushner in the 1970s to approximate optimally controlled diffusion processes \cite{Ku1970, Ku1973, Ku1975, Ku1976A, Ku1976B, Ku1977, Ku2001, KuDu2001, Ku2002,Ku2004, Ku2006, Ku2011}.  However, because of their interest in stochastic control problems, these works focus on numerical solutions with gridded state spaces  that admit a global matrix representation.  In the statistical physics literature, this matrix representation was avoided and a Monte-Carlo method was used to simulate the numerical solution \cite{elston1996numerical, WaPeEl2003,Ph2008,MeScVa2009,LaMeHaSc2011}, and to be specific, this idea seems to go back to at least \cite{elston1996numerical}.  Among these papers, the most general and geometrically flexible MCAM were the finite volume methods developed for over-damped Langevin equations presented in \cite{LaMeHaSc2011}.  More recently, the MCAM framework has been generalized to lessen the requirements on the underlying diffusion process \cite{BoVa2016}.  In particular,  this generalization no longer requires that the domain of the diffusion process is bounded, that the infinitesimal generator of the diffusion process is symmetric or that the infinitesimal covariance of the diffusion process is diagonally dominant.  This generalization is made possible by letting the state space of the numerical solution be gridless and by using Monte-Carlo methods to simulate the numerical solution, but at the same time, keeping the restriction that the jump size of the approximation is uniformly bounded.     As we will see below, this property is essential to the stability and accuracy of SPECTRWM.

\section{Algorithm} \label{sec:spectrwm}

We present two versions of SPECTRWM.  The point of the first version is to illustrate basic concepts, and as such, it uses a spatial finite difference approximation of an SPDE on an evenly spaced grid in 1D. (This `academic' version was used by the author in his mini-course entitled ``Spectral Random Walk Method for the Numerical Solution of Stochastic Partial Differential Equations'' at the 2016 Gene Golub SIAM Summer School.)   The second version is based on a spectral Galerkin approximation, which, as we will see, is more general and more efficient than the academic one.
Unless otherwise stated, we will mainly provide numerical verification of the academic version of SPECTRWM.

\subsection{Academic Version}

We present this version of SPECTRWM in the specific context of a one-dimensional SPDE with additive, space-time Gaussian white noise and a scalar noise coefficient.   We assume that we are given a stable and accurate semi-discrete system that consists of a 1D grid with $n$ grid points, an initial condition on this grid, a spatial step size parameter $\Delta x>0$, and a system of  approximating SDEs on $\mathbb{R}^n$ of the form:
\begin{equation} \label{eq:semi_discrete}
d u^n =  \left[ L_n u^n + F_n ( u^n )  \right] dt + \sigma  \sqrt{\frac{1}{\Delta x}} dW_n
\end{equation}
where $L_n$ is an $n \times n$ discretization matrix associated to the linear part of the drift, $F_n$ is a discretization of the nonlinear part of the drift, $\sigma>0$ is a scalar constant, and $W_n$ is an $n$-dimensional Brownian motion.  We assume that $L_n$ has an orthonormal set of $n$ eigenvectors $\{ e_i \}$ with eigenvalues $\{ \mu_i \}$.   Implicit in this setup are the boundary conditions of the SPDE, which are typically built into \eqref{eq:semi_discrete}.  A typical example of $L_n$ is the discrete Laplacian for the standard finite difference method with periodic, Dirichlet, Neumann or mixed boundary conditions.  In what follows, $u^h(t) \in \mathbb{R}^n$ denotes the Markov jump process approximation produced by SPECTRWM.  

\medskip

\begin{algorithm}[SPECTRWM: Academic Version] \label{algo:spectrwm_academic}
 Given the current time $t \ge 0$, the current state of the process $u_h(t)=v \in \mathbb{R}^n$, a spatial step size $\Delta x>0$, and a jump size $h>0$,  the algorithm outputs an updated state  $u_h(t+\delta t)$ at time $t+\delta t$ in three sub-steps.
 \begin{description}
\item[(Step 1)] compute forward/backward jump rates: \begin{equation}  \label{eq:jump_rates_academic}
J^{\pm}_i(v) = \frac{\sigma^2}{2 h^2 \Delta x} \exp\left( \pm \left( \mu_i v^T e_i + F_n(v)^T e_i \right) \frac{h \Delta x}{\sigma^2} \right) 
\end{equation} for $1 \le i \le n$.  
\item[(Step 2)] update time via \[
t_1 = t_0 + \delta t
\] where $\delta t$ is an exponentially distributed random variable with parameter \begin{equation} \label{eq:sumJrates}
J(v) = \sum_{i=1}^n ( J_i^+(v) + J_i^-(v) ) \;.
\end{equation}
\item[(Step 3)] update the state of the system by assuming that the process jumps forward/backward along the eigenvector $e_i$ to state $v \pm h e_i$ with probability: \[
\Pr( u_h(t + \delta t) = v \pm h e_i  \mid u_h(t)=v )  = \frac{J_i^{\pm}(v)}{J(v)} 
\]  
for $1 \le i \le n$. 
\end{description}
\end{algorithm}

\medskip

We stress that this method is very simple and straightforward to implement.  Note from (Step 3) that the algorithm moves by jumps in the semi-discrete space $\mathbb{R}^n$ in the directions of the eigenvectors of $L$.  Moreover, the jump size is $h$, and these jumps alter every component, i.e., each jump induces a global system update.    The time elapsed in each state $\delta t$ is an exponentially distributed random variable with parameter $J(v)$ given in \eqref{eq:sumJrates}, which is defined as the sum of the forward/backward jump rates from (Step 1).  
Also, the update rules in (Step 1) and (Step 2) only depend on the current state.  Thus, the resulting process is a Markov jump process \cite{ReYo1999,protter2004stochastic,EtKu2009,Kl2012}.

This process has an infinitesimal generator that is given by: \begin{equation} \label{eq:generator}
Qf(v) = \sum_{1 \le i \le n} J_i^+(v) \left( f(v+h e_i) - f(v) \right) + J_i^-(v) \left( f(v- h e_i) - f(v) \right) \;.
\end{equation}
Using a Taylor expansion of $Qf(v)$ about $h=0$, in Appendix~\ref{sec:local_consistency} we show that: \begin{equation} \label{eq:local_consistency}
Qf(v)=  v^T L_n \nabla f(v) + F_n(v)^T \nabla f(v) + \frac{\sigma^2}{2 \Delta x} \tr( D^2 f(v) ) + O(h^2 \Delta x + h^2 ) \;.
\end{equation} We recognize the leading order term in this expansion as the infinitesimal generator of the approximating SDE in \eqref{eq:semi_discrete}.  Thus, by choosing $h$ sufficiently small we can get arbitrarily close (in law) to the solution of the approximating SDE.

However, a very practical question remains: what is the computational cost of SPECTRWM?  To address this question, we consider the scaling of SPECTRWM as the jump size decreases $h \to 0$; and, as the spatial step size decreases $\Delta x \to 0$ (or equivalently $n \to \infty$).  In order to produce a single trajectory over $[0, t]$, the computational cost of SPECTRWM is proportional to the (random) number of steps $N(t)$ it takes for the approximation to reach $t$.  In \S4.3 of \cite{BoVa2016}, it is shown that the mean of $N(t)$ is inversely proportional to the mean holding time.  Thus, the mean holding time dictates the average cost of the algorithm, which from (Step 2) of Algorithm~\ref{algo:spectrwm} scales like $O(h^2 / n^2 )$.   This scaling of the mean holding time reflects the roughness of the noise and the fact that SPECTRWM is able to jump in the direction of any of the $n$ eigenmodes of $L$.   However, it turns out we can do much better than this, by using a spectral Galerkin approximation, as we describe next.


\subsection{Fast Version}

Consider an SPDE with solution $u$ defined on a Hilbert space $H$ with inner product $\langle \cdot , \cdot \rangle$.  Assume that the noise in the SPDE is additive space-time Gaussian white noise with scalar noise coefficient $\sigma>0$.   Turning to the drift of the SPDE, assume this drift has a linear part of the form $Lu$ where $L$ is a linear operator with a complete orthonormal set of eigenfunctions $\{ e_i \}$ and eigenvalues $\{ \mu_i \}$ ordered such that $\mu_i \ge \mu_{i+1}$ for all natural numbers $i  > 0$.  Typically $L$ is a uniformly elliptic differential operator with non-positive eigenvalues.  In terms of these eigenfunctions, we define the finite-dimensional subspace $V_n = \spanset\{ e_1, \cdots, e_n \}$, which is the span of the eigenfunctions associated to the $n$ largest eigenvalues of $L$, and let $u^n$ denote the orthogonal projection of $u$ onto this finite-dimensional subspace $V_n$: \[
u^n = \sum_{i=1}^n \hat{u}_i^n e_i  \;.
\] where we introduced the $n$ spectral coefficients: $\{ \hat{u}_i^n = \langle u^n, e_i \rangle~\mid~ 1 \le i \le n \}$.   The function $u^n$ is a standard spectral Galerkin approximation of $u$.  Next we derive an $n$-dimensional system of approximating SDEs that the spectral coefficients of $u^n$ satisfy.

For this purpose, we approximate the space-time Wiener process in the SPDE by a truncated sum: $\sum_{i=1}^n B_i(t) e_i $
 where $\{ B_i \}$ are iid Brownian motions.  Finally, let $\{ \hat{F}_i \}$ denote the spectral coefficients of the rest of the drift of the SPDE in the finite-dimensional basis $\{ e_i \}_{1 \le i \le n}$.   With these truncations in hand, the spectral coefficients of $u^n$ satisfy the following system of SDEs: \begin{equation} \label{eq:spectral_galerkin_approx}
d \hat{u}_i^n = \left[  \mu_i \hat{u}_i^n dt + \hat{F}_i( u^n ) \right] dt + \sigma dB_i 
\end{equation}
where $i$ ranges from $1$ to $n$.    The fast version of SPECTRWM directly approximates these SDEs in the spectral domain.  Let $u^h(t)$ denote the numerical solution produced by SPECTRWM.

\medskip

\begin{algorithm}[SPECTRWM: Fast Version] \label{algo:spectrwm}
 Given a jump size $h>0$, the current time $t \ge 0$, and the current state of the process $u^h(t) = v$,  the algorithm outputs an updated state $u^h(t+\delta t)$ at time $t+\delta t$ in three sub-steps.  
 \begin{description}
\item[(Step 1)] compute forward/backward jump rates: \begin{equation}  \label{eq:jump_rates_fast}
J^{\pm}_i( v )  = \frac{\sigma^2}{2 h^2} \exp\left( \pm \left(  \mu_i \hat{v}_i + \hat{F}_i( v )  \right) \frac{h}{\sigma^2} \right) 
\end{equation} for $1 \le i \le n$.  
\item[(Step 2)] update time via \[
t_1 = t_0 + \delta t
\] where $\delta t$ is an exponentially distributed random variable with parameter \begin{equation} \label{eq:sumJrates}
J( v ) = \sum_{i=1}^n \left[ J_i^+( v ) + J_i^-( v ) \right] \;.
\end{equation}
\item[(Step 3)] update the state of the system by assuming that the process jumps forward/backward along the eigenvector $e_i$ with probability: \[
\Pr( u^h(t+\delta t) = v \pm h e_i   \mid u^h(t) = v )  = \frac{J_i^{\pm}( v )}{J( v )} 
\]  
for $1 \le i \le n$. 
\end{description}
\end{algorithm}

\medskip

Note that this version of SPECTRWM makes jumps of size $h$ in the individual spectral coefficients of $u^h$.  Every jump in the spectral domain leads to a jump of size $h$ in $u^h$ along the corresponding eigenfunction.  The generator of this version of SPECTRWM is identical in form to the generator given in \eqref{eq:generator} with forward/backward jump rates given in \eqref{eq:jump_rates_fast}.  However, in this case, the jumps and state of the system are in the spectral domain where the approximating SDEs are defined.  In this case, a straightforward Taylor expansion, shows that the infinitesimal generator of this version is an $O(h^2)$ approximation to the infinitesimal generator associated to the approximating SDEs in \eqref{eq:spectral_galerkin_approx}.    This version of SPECTRWM is `fast' because the mean holding time in this approximation scales like $O(h^2 / n)$, in contrast to the $O(h^2 / n^2)$ scaling of the academic version.

\section{Heat SPDE}

To assess accuracy of SPECTRWM, consider the heat SPDE on $[0, 2 \pi]$ with periodic boundary conditions at  $0$ and $1$: \begin{equation} \label{eq:heat_spde}
\begin{dcases}
du =  \left( \dfrac{\partial^2 u}{\partial x^2}  -\lambda u \right) dt + \sigma d W & \forall x \in [0,2 \pi], t \ge 0 \\
u(t,0) =  u(t,2 \pi)  & \forall t\ge 0 \;, \\
u(0,x) = \frac{c_0}{\sqrt{2 \pi}} + \frac{1}{\sqrt{\pi}} \sum_{k \ge 1} \left( c_{-k} \cos(k x) + c_k \sin(k x) \right)  & \forall x \in [0,2 \pi] \;.
\end{dcases}
\end{equation}
where $\lambda \ge 0$ and $\sigma>0$ are parameters; and $\{ c_k \}$ are the Fourier coefficients of the initial conditions.  By the usual Fourier series argument, the solution to these equations at any time $t>0$ is a Gaussian process with mean \begin{equation} \label{eq:mean_heat_spde}
\E u(t,x)  =  \frac{c_0}{\sqrt{2 \pi}} e^{-\mu_0 t} + \frac{1}{\sqrt{\pi}} \sum_{k \ge 1} e^{- \mu_k t} \left( c_{-k} \cos(k x) + c_k \sin(k x) \right) 
\end{equation} and spatial covariance \begin{equation} \label{eq:spatial_covariance_heat_spde}
\begin{aligned}
& \E\left\{ (u(t,x) - \E u(t,x))  (u(t,y) - \E u(t,y))  \right\} = \\
& \qquad\qquad  \begin{dcases}
 \frac{\sigma^2}{2 \pi} t + \frac{\sigma^2}{\pi} \sum_{k \ge 1} \frac{1 - e^{-2 k^2 t}}{2 k^2} \cos(k (x-y) ) & \text{if $\lambda=0$} \\
  \frac{\sigma^2}{2 \pi} \frac{1 - e^{-2 \mu_0 t}}{2 \mu_0} + \frac{\sigma^2}{\pi} \sum_{k \ge 1} \frac{1 - e^{-2 \mu_k t}}{2 \mu_k} \cos(k (x-y) ) & \text{if $\lambda>0$} 
 \end{dcases}
 \end{aligned}
\end{equation} 
where for convenience we have introduced $\mu_k = k^2 + \lambda$.  
We use a truncation of \eqref{eq:mean_heat_spde} and \eqref{eq:spatial_covariance_heat_spde} as a benchmark solution to verify the efficiency of the method in approximating expected values involving (a) the solution at a fixed time and (b) the entire path of the solution.  

\begin{remark}
The dynamics of the heat SPDE involves a `competition' between the dissipative drift term and the space-time Gaussian white noise.  When $\lambda=0$, the noise wins in the sense that the zeroth order mode of the Fourier series solution is a Brownian motion, and as a consequence, the SPDE has no stationary distribution.  This is reflected in the secular term appearing in \eqref{eq:spatial_covariance_heat_spde} in the case $\lambda=0$.   However, for any $\lambda>0$, the zeroth order mode of the Fourier series solution is an OU process, and the heat SPDE has a stationary distribution.  
\end{remark}

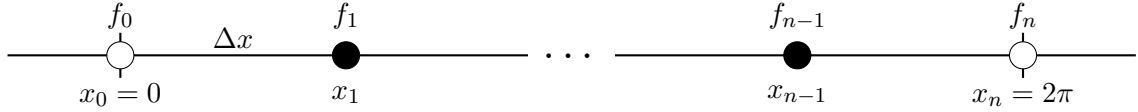
\begin{figure}
\begin{tikzpicture}[scale=1.5]
\draw[-, thick](0,0.0) -- (10,0.0);
\draw[-, thick](1,-0.2) -- (1,0.2);
\draw[-, thick](9,-0.2) -- (9,0.2);
\node[black,scale=1] at (2,0.15) {$\Delta x$};
\node[black,scale=1] at (1,-0.35) {$x_0=0$};
\node[black,scale=1] at (9,-0.35) {$x_n=2 \pi$};
\node[black,scale=1.] at (3,-0.35) {$x_1$};
\node[black,scale=1.] at (1,0.35) {$f_0$};
\node[black,scale=1.] at (3,0.35) {$f_1$};
\node[black,scale=1.] at (7,0.35) {$f_{n-1}$};
\node[black,scale=1.] at (9,0.35) {$f_{n}$};
\node[black,scale=1.] at (7,-0.35) {$x_{n-1}$};
\node[black, scale=1.5,fill=white] at (5.0,0.0) {$\dotsc$};
\filldraw[color=black,fill=white] (1,0) circle (0.12);
\filldraw[color=black,fill=black] (3,0) circle (0.12);
\filldraw[color=black,fill=black] (7,0) circle (0.12);
\filldraw[color=black,fill=white] (9,0) circle (0.12);
\end{tikzpicture}
\caption{\small  {\bf Grid.}  This cartoon shows the grid used for all of the test problems in this paper.  The black dots mark the interior grid points, and the white dots mark the boundary grid points.  Due to the periodic boundary conditions, the approximating SDEs only need to be evolved on the $n$ grid points: $x_0, \cdots , x_{n-1}$.  Also, for any function $f: [0, 2\pi] \to \mathbb{R}$, we use the shorthand notation $f_i = f(x_i)$ for $0 \le i \le n$.  
}
  \label{fig:grid}
\end{figure}

\begin{figure}[ht!]
\begin{center}
\includegraphics[width=0.45\textwidth]{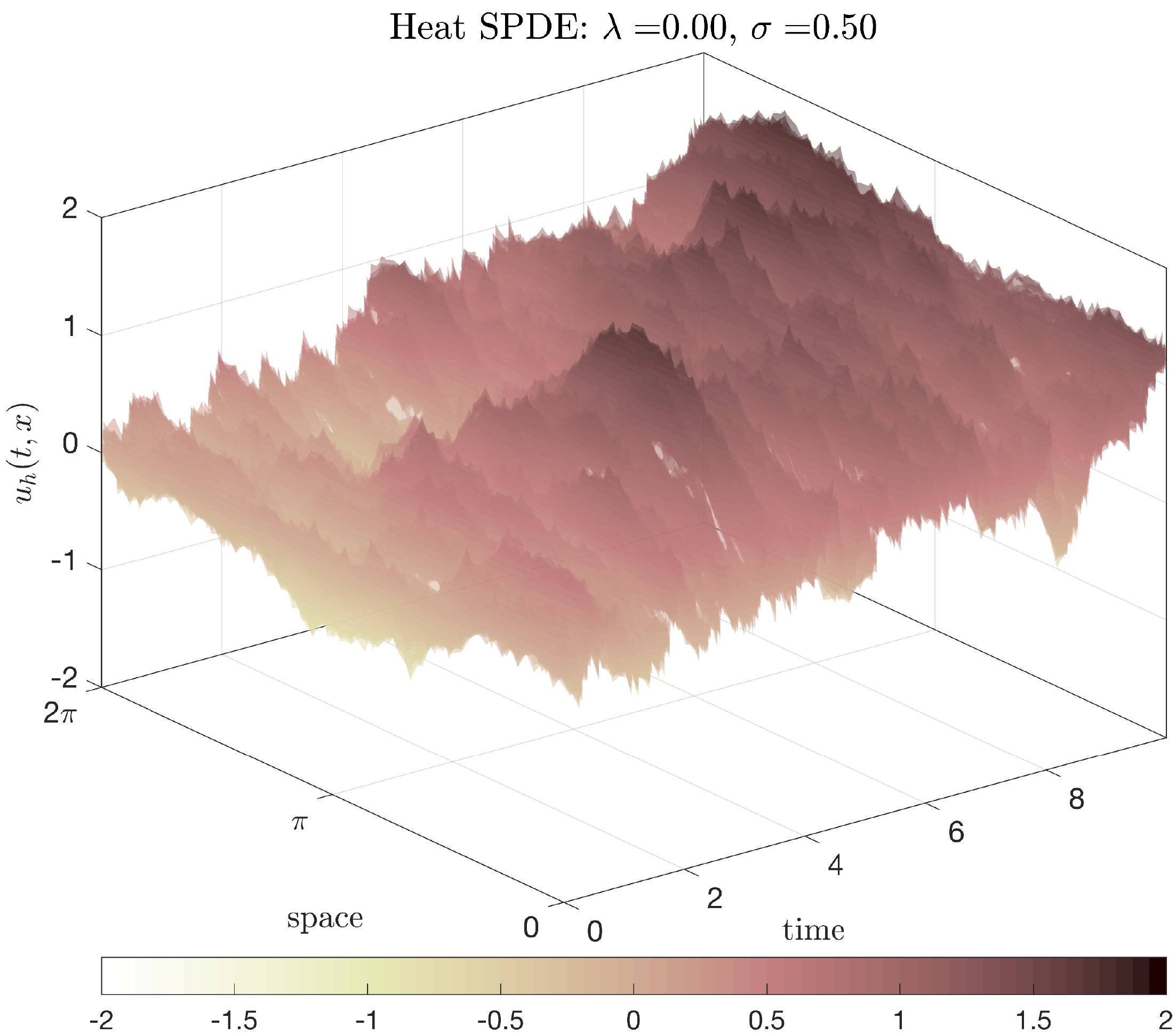} 
\includegraphics[width=0.45\textwidth]{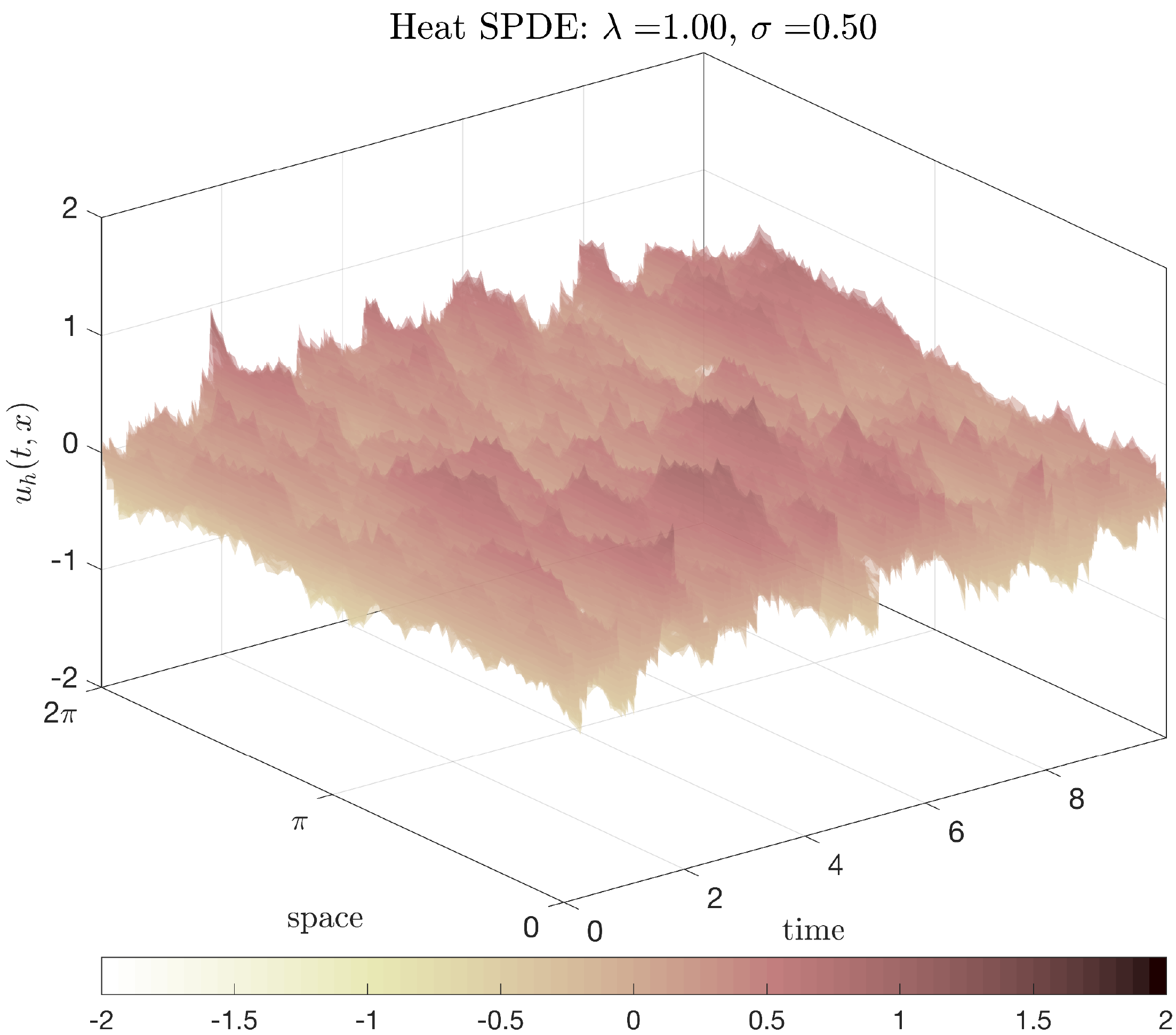} 
\end{center}
\caption{ \small {\bf Heat SPDE on $[0,2 \pi]$.}   Sample paths produced using SPECTRWM at the parameter values indicated in the figure title and with a trivial initial condition.  Colors indicate the height of the surface, and are added in order to make the surface plots a bit more clear.  These sample paths illustrate the stabilizing effect of the parameter $\lambda$.   }
 \label{fig:spectrwm_heat_spde_sample_paths}
\end{figure}

\begin{figure}[ht!]
\begin{center}
\includegraphics[width=0.45\textwidth]{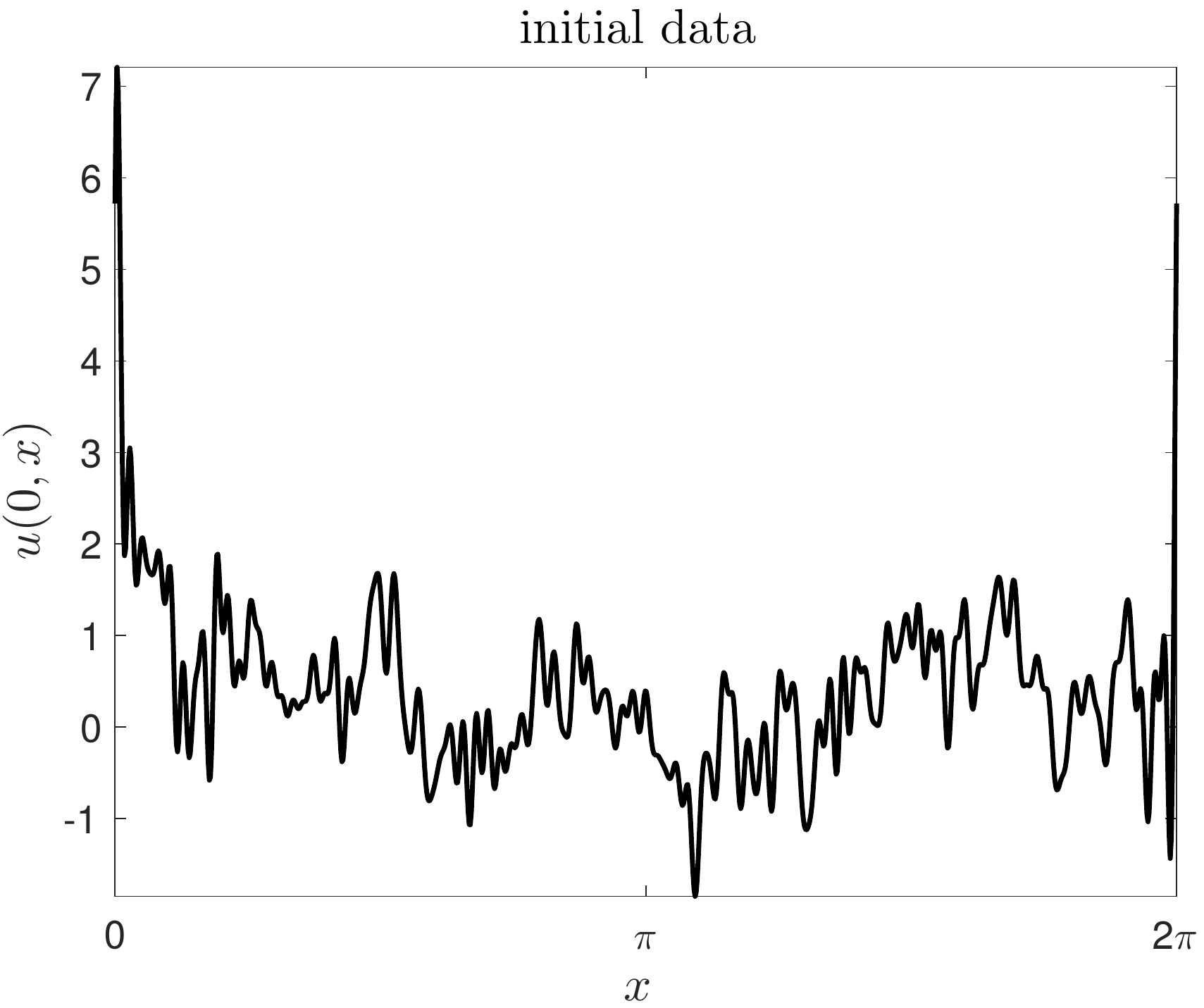} 
\includegraphics[width=0.45\textwidth]{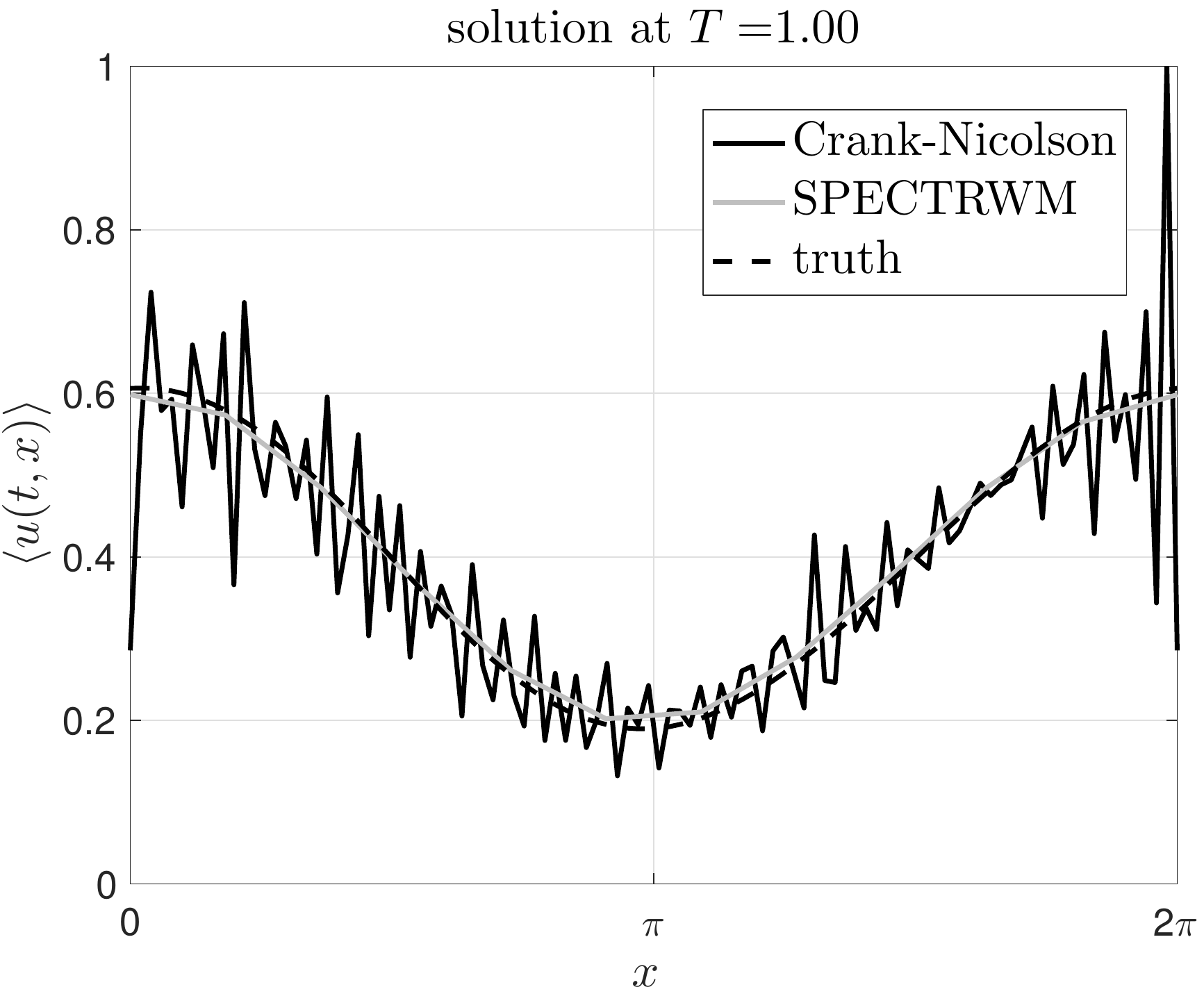} 
\end{center}
\caption{ \small {\bf Heat SPDE on $[0,2 \pi]$.}
This figure compares the fast version of SPECTRWM with $\sigma=1$, $\lambda=0$, and $n=11$ with the Crank-Nicolson scheme using $n=101$ grid points and $\Delta t = \Delta x$.  The initial condition for this numerical test is shown in the left panel.  The expected value of the true solution is the dashed line shown in the right panel.  As is well known, the Crank-Nicolson scheme has trouble damping the high frequency modes, whereas the fast version of SPECTRWM does not have this issue.
}
 \label{fig:spectrwm_comparison}
\end{figure}

Referring to Figure~\ref{fig:grid}, we discretize $[0,2 \pi]$ using an evenly-spaced grid: \[
\left\{ x_i =  i \Delta x   \mid 0 \le i \le n \right\}
\]  with spatial step size $\Delta x =2 \pi /n$.  Let $f_i=f(x_i)$.  On this grid, we approximate the diffusive part of the drift in \eqref{eq:heat_spde} by a finite-difference method: \[
\frac{\partial^2 f}{\partial x^2}(x_i) = \dfrac{f_{(i+1) \bmod n} - 2 f_i + f_{(i-1) \bmod n}}{\Delta x^2} + O(\Delta x^2)
\]
which is valid if $f$ has four derivatives.  The resulting discretization matrix is the discrete Laplacian with periodic boundary conditions in 1D, and its eigenvalues and eigenvectors are analytically known.  Samples paths generated by SPECTRWM are given in Figure~\ref{fig:spectrwm_heat_spde_sample_paths} for two different values of $\lambda$.  The finite time weak accuracy of the scheme is verified in the graphs labelled (a) in the legend of Figure~\ref{fig:spectrwm_accuracy}.  The accuracy of the method in approximating the time integral of this covariance over $[0, 1]$ (a path-dependent expected value) is verified in the graph labelled (b) in the figure legend.   Both of these graphs are in agreement with the local consistency estimate provided in \eqref{eq:local_consistency}.   Figure~\ref{fig:spectrwm_comparison} shows that SPECTRWM is able to accurately represent the decay of high frequency modes.  This is expected since the algorithm is based on a random walk in the spectral coefficients.


\begin{figure}[ht!]
\begin{center}
\includegraphics[width=0.45\textwidth]{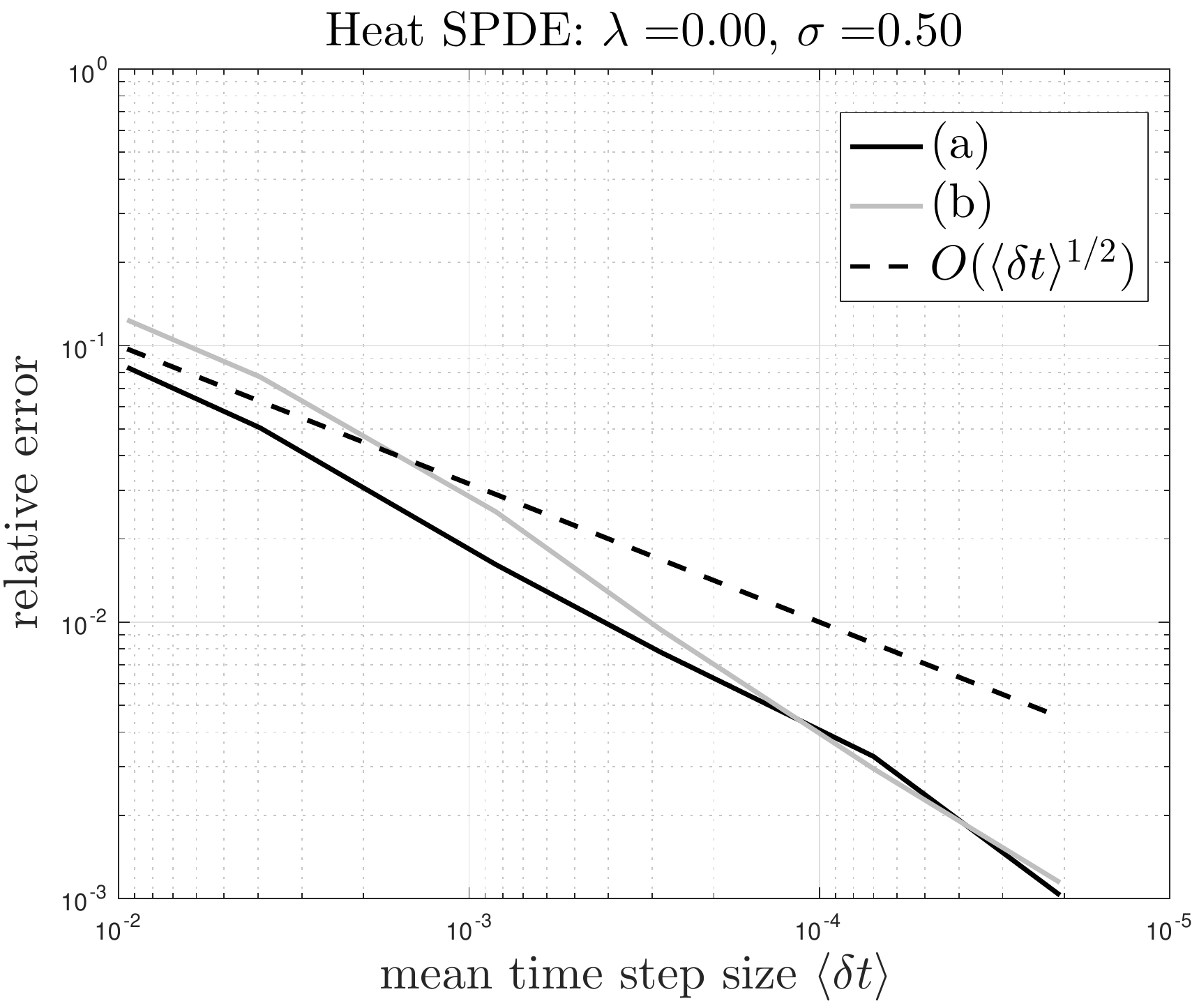} 
\includegraphics[width=0.45\textwidth]{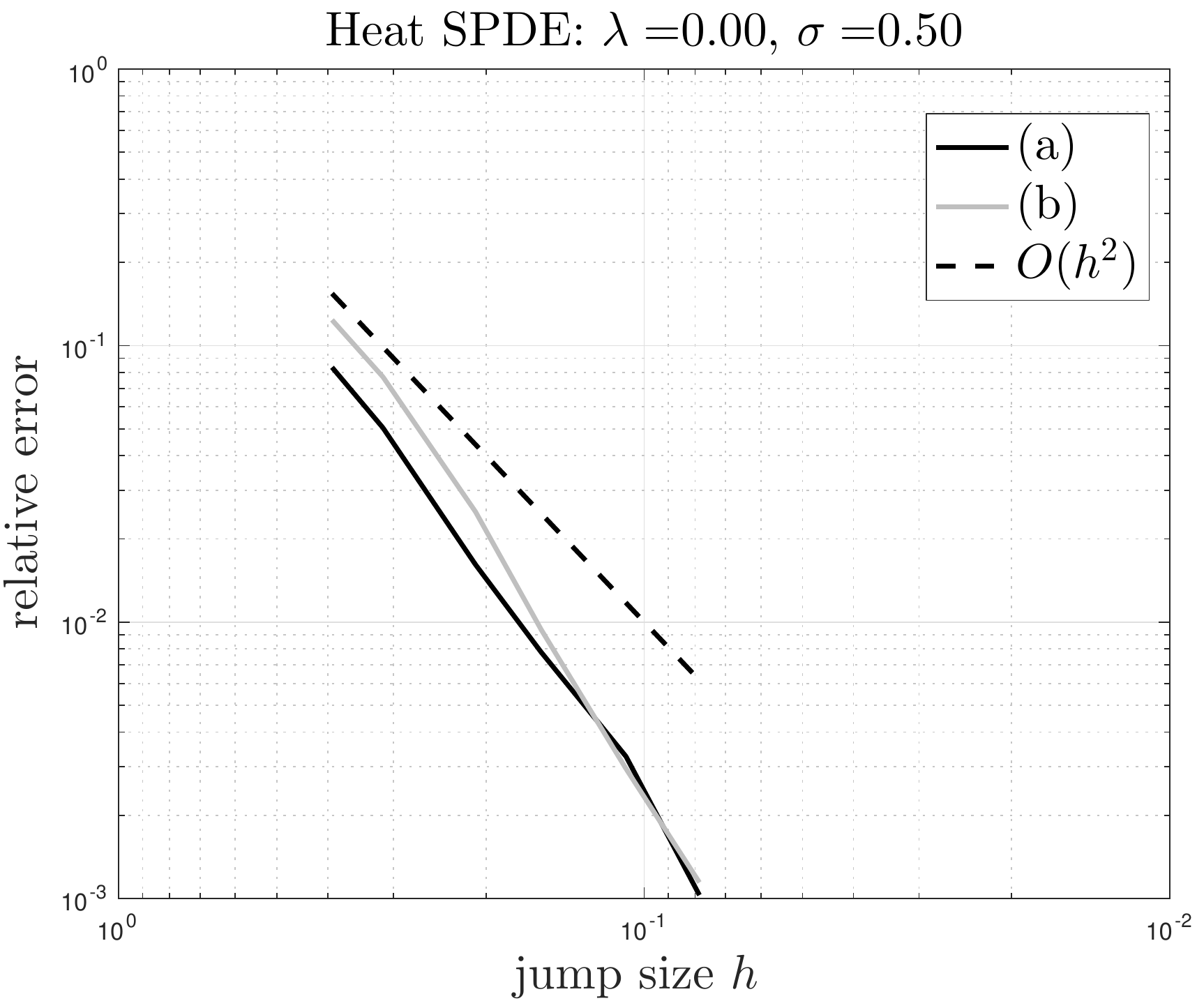} 
\end{center}
\caption{ \small {\bf Heat SPDE on $[0,2 \pi]$.}
This figure verifies the accuracy of SPECTRWM for the parameter values given in the figure title and for (a) the expected value of the solution squared at time $t=1$ and (b) the time integral of this expected value over $[0,1]$.  The dashed lines are for reference.  
}
 \label{fig:spectrwm_accuracy}
\end{figure}


\section{Overdamped Langevin SPDE}

\begin{figure}[ht!]
\begin{center}
\includegraphics[width=0.75\textwidth]{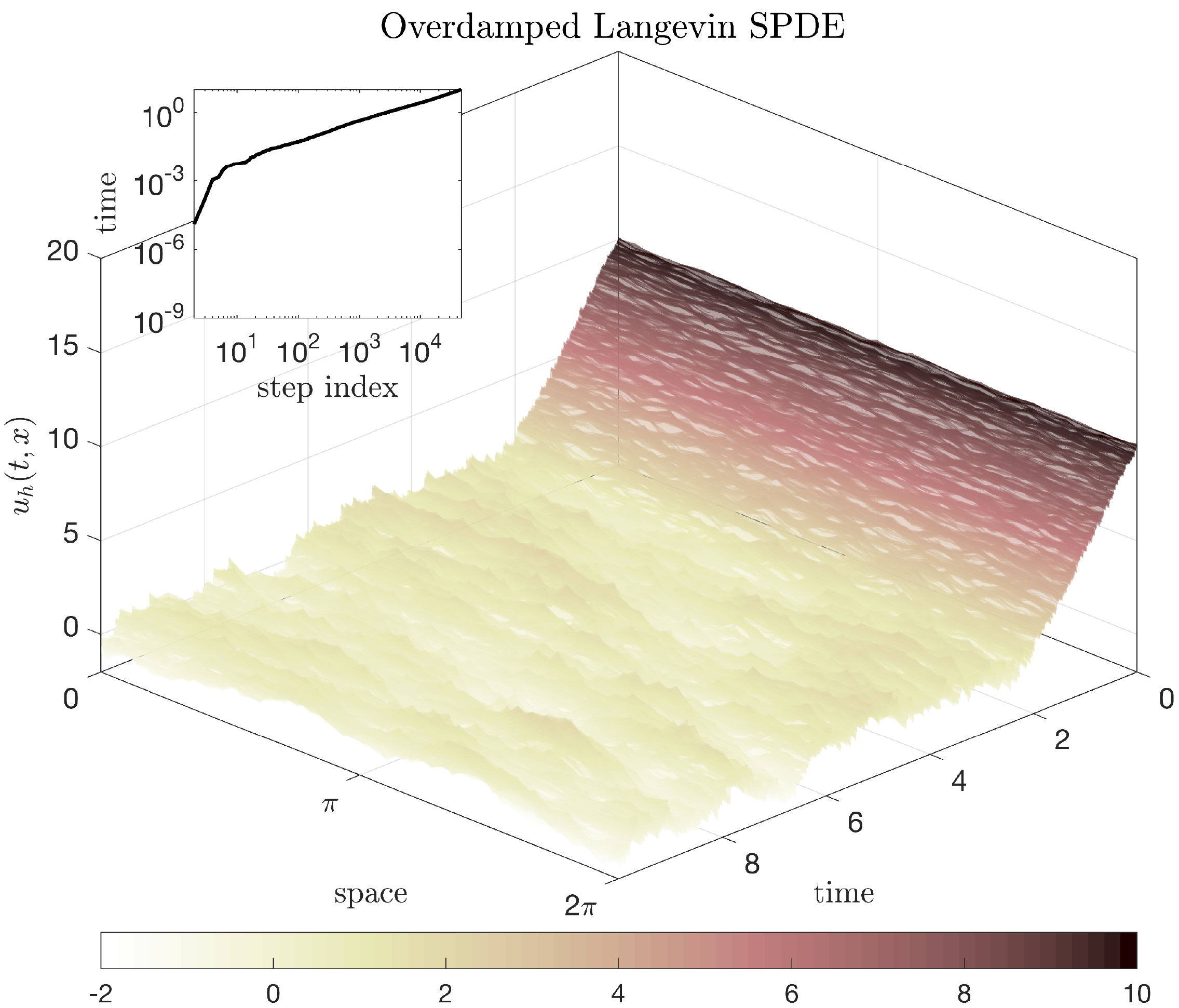}  
\end{center}
\caption{ \small {\bf Overdamped Langevin SPDE on $[0,2 \pi]$.}   Sample paths produced using SPECTRWM at the parameter values indicated in the figure title.  Colors indicate the height of the surface, and are added in order to make the surface plots a bit more clear.  These sample paths illustrate the stabilizing effect of the parameter $\lambda$.   }
 \label{fig:spectrwm_langevin_spde_sample_paths}
\end{figure}

To assess ergodicity of SPECTRWM, consider a nonlinear, overdamped Langevin SPDE: \begin{equation}
\label{eq:heat_spde}
\begin{dcases}
du =  \left( \dfrac{\partial^2 u}{\partial x^2}  - u^3 \right) dt + \sigma d W & \forall x \in [0,2 \pi], t \ge 0 \\
u(t,0) =  u(t,2 \pi)  & \forall t\ge 0 \;, \\
\end{dcases}
\end{equation}
with initial conditions that we will describe shortly.  Using the same discretization as before, the approximating SDEs are themselves overdamped Langevin SDEs, which preserve a probability density function: \[
\pi^n(v) \propto \exp\left( -\Delta x \left( \sum_{0 \le i \le n-1} \frac{v_i^4}{4} + \frac{v^T L v}{2} \right) \right) 
\] where $v = (v_0, \cdots, v_{n-1})$.  For the first test, we took an initial condition at very high energy as shown in Figure~\ref{fig:spectrwm_langevin_spde_sample_paths}.   This figure suggests that SPECTRWM is geometrically ergodic when the underlying SDE is.  For the second test, we took a trivial initial condition, a very large jump size $h=\sqrt{\Delta x}$ and computed the accuracy of SPECTRWM. Following Chapter 2 of \cite{BoVa2016}, we modified the jump rates so that the algorithm exactly preserves the stationary density of the approximating SDEs.

\begin{figure}[ht!]
\begin{center}
\includegraphics[width=0.45\textwidth]{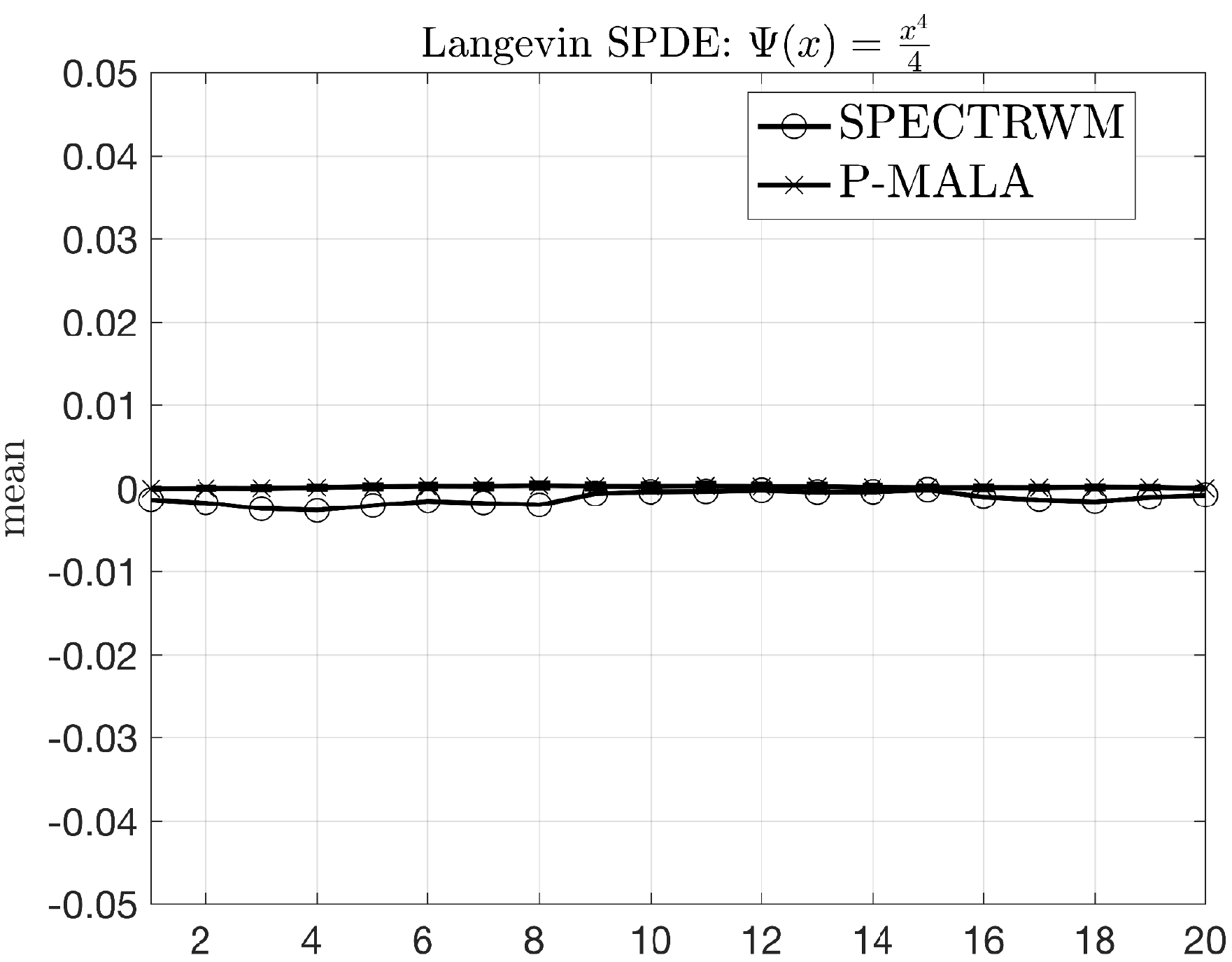} 
\includegraphics[width=0.45\textwidth]{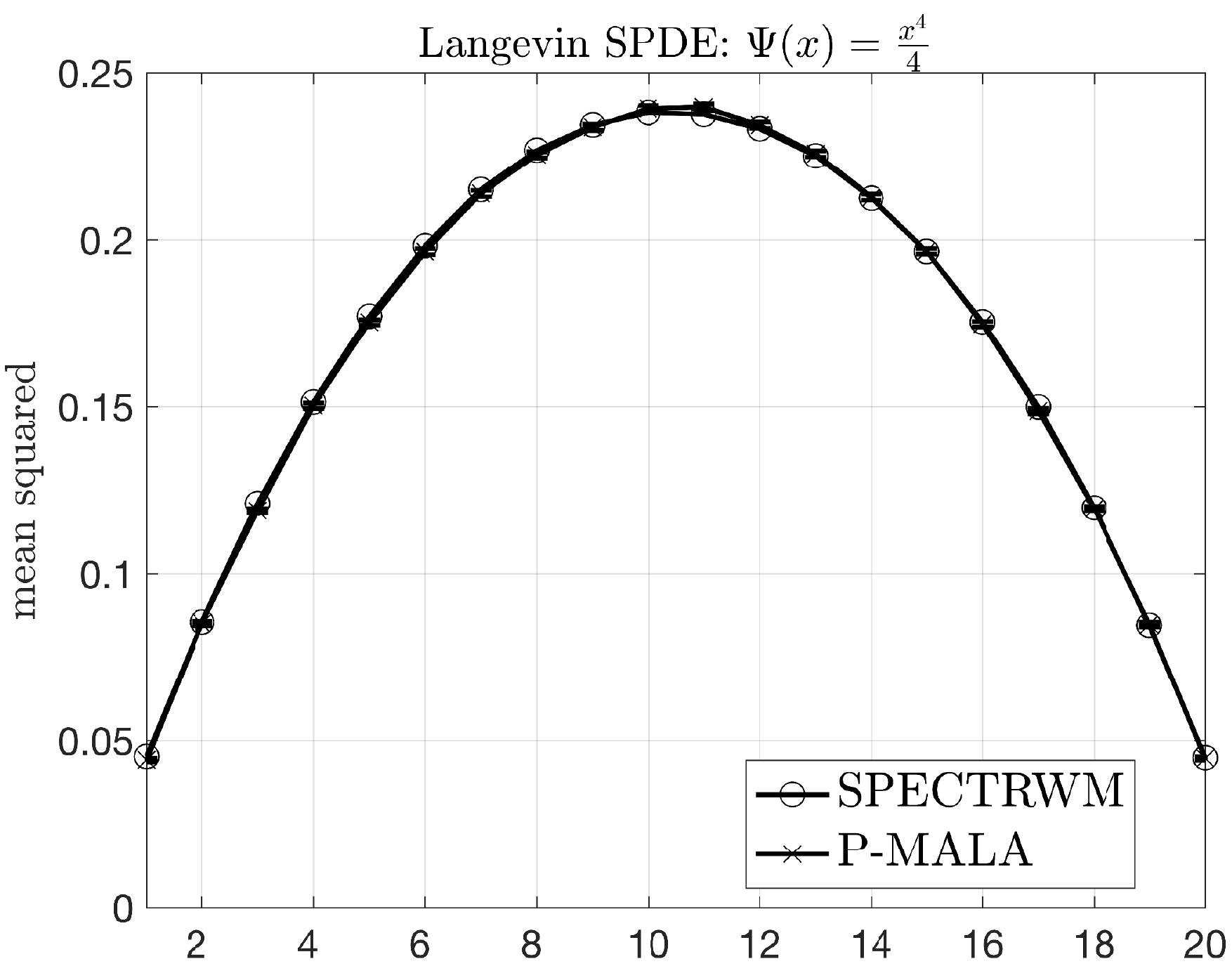} 
\end{center}
\caption{ \small {\bf Overdamped Langevin SPDE on $[0,2 \pi]$.}
This figure verifies the accuracy of SPECTRWM in computing the first and second moment of each component of the approximating SDE with $N=20$ grid points.  The $x$-axis labels the components.  We use as benchmark solution the preconditioned MALA algorithm, which obtains proposal moves from the $\theta$-method with $\theta=1/2$ \cite{BeRoStVo2008}.
}
 \label{fig:spectrwm_accuracy}
\end{figure}

\section{Burgers SPDE}

\begin{figure}[ht!]
\begin{center}
\includegraphics[width=0.65\textwidth]{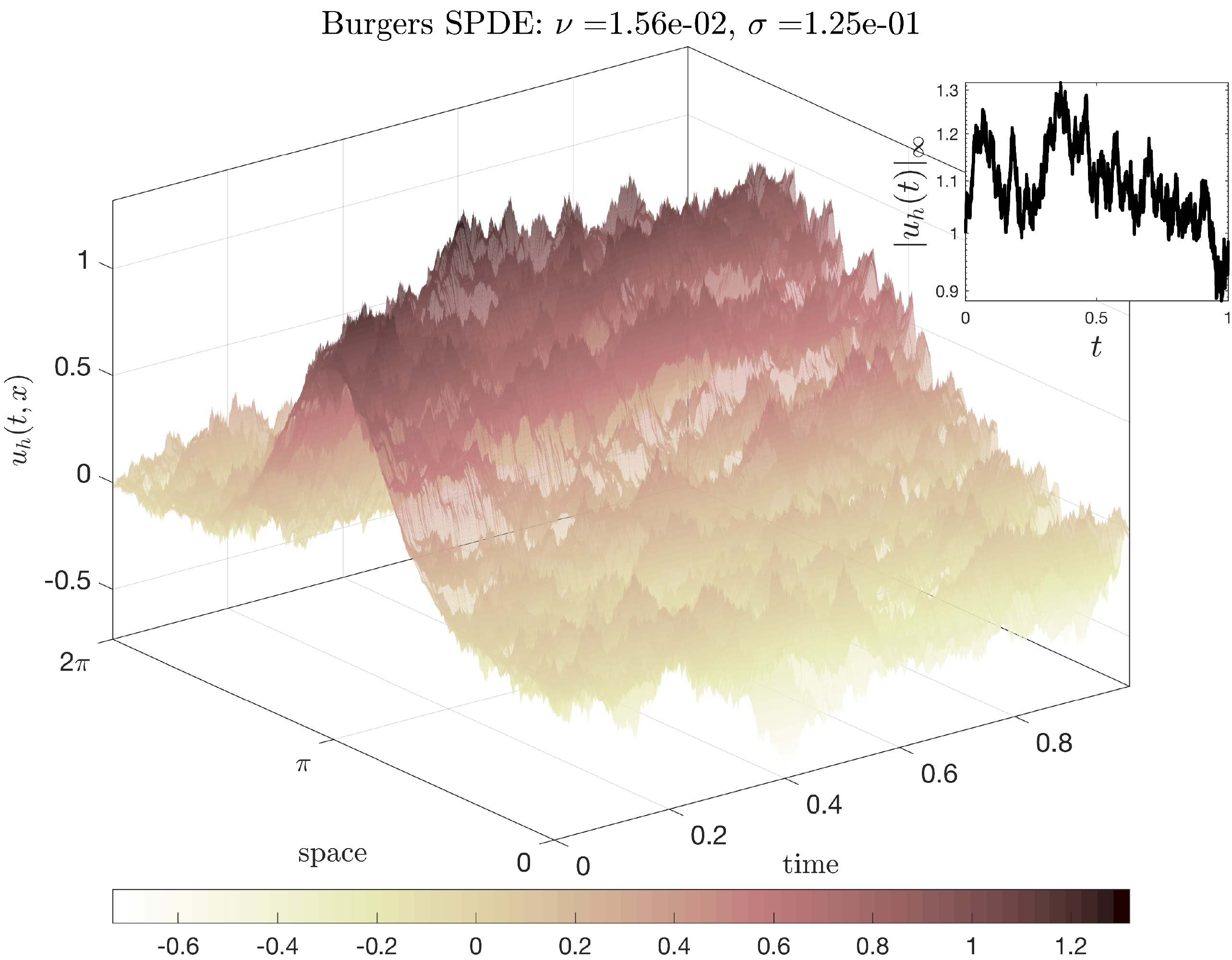}  \\
\includegraphics[width=0.65\textwidth]{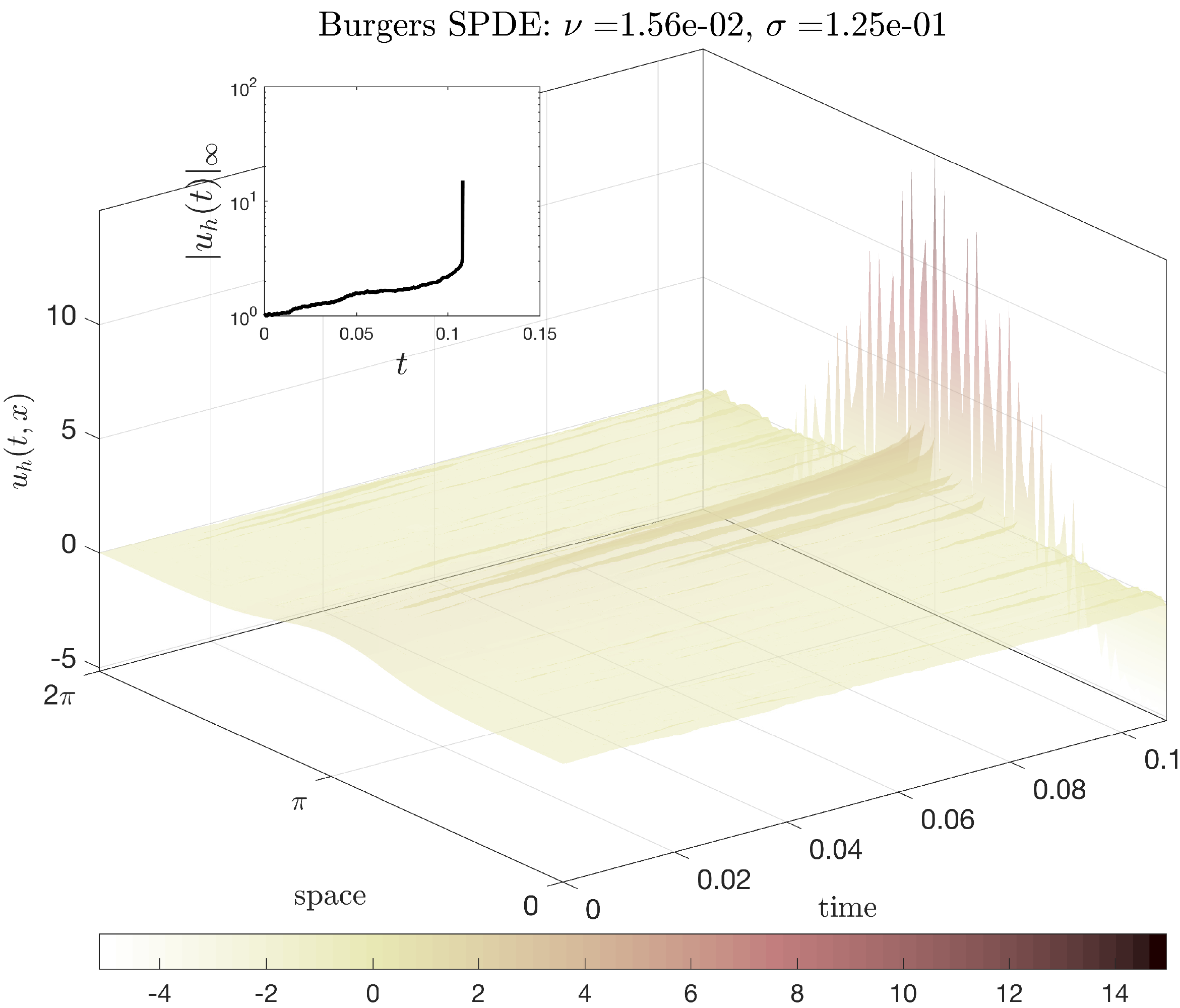} 
\end{center}
\caption{ \small {\bf Burgers SPDE on $[0,2 \pi]$.}   Sample paths produced using SPECTRWM at the parameter values indicated in the figure title.  Colors indicate the height of the surface, and are added in order to make the surface plots a bit more clear.  The left panel uses a valid semi-discrete approximation, while the right panel uses a semi-discrete approximation that converges to a non-ergodic SPDE.  }
 \label{fig:spectrwm_burgers_spde_sample_paths}
\end{figure}

We apply SPECTRWM to the following Burgers SPDE: \begin{equation}
\label{eq:heat_spde}
\begin{dcases}
du =  \left( \nu \dfrac{\partial^2 u}{\partial x^2}  - u\dfrac{\partial u}{\partial x}   \right) dt + \sigma d W & \forall x \in [0,2 \pi], t \ge 0 \\
u(t,0) =  u(t,2 \pi)  & \forall t\ge 0 \;,
\end{dcases}
\end{equation}
with a bump function initial condition as depicted in Figure~\ref{fig:spectrwm_burgers_spde_sample_paths}.
To construct approximating SDEs, we use the same setup for the linear part and discretize the nonlinear advective term using either:\[
f \dfrac{\partial f}{\partial x}(x_i)  \approx \begin{dcases}
\frac{1}{2 \Delta x} \left( f^2_{(i+1) \bmod n}- f^2_{(i-1) \bmod n} \right) \\
\frac{1}{\Delta x}  \left( f_{(i+1) \bmod n} - f_{i} \right) f_i
\end{dcases}
\]  As we confirm in Figure~\ref{fig:spectrwm_burgers_spde_sample_paths}, the latter discretization leads to an approximation that is unstable.  This qualitatively confirms the results in \cite{HaVo2011, HaMa2012}.  A quantitative comparison is not possible, because the precise form of the spurious drift term depends on the approximation being used.  

\section{KPZ SPDE} 

\begin{figure}[ht!]
\begin{center}
\includegraphics[width=0.65\textwidth]{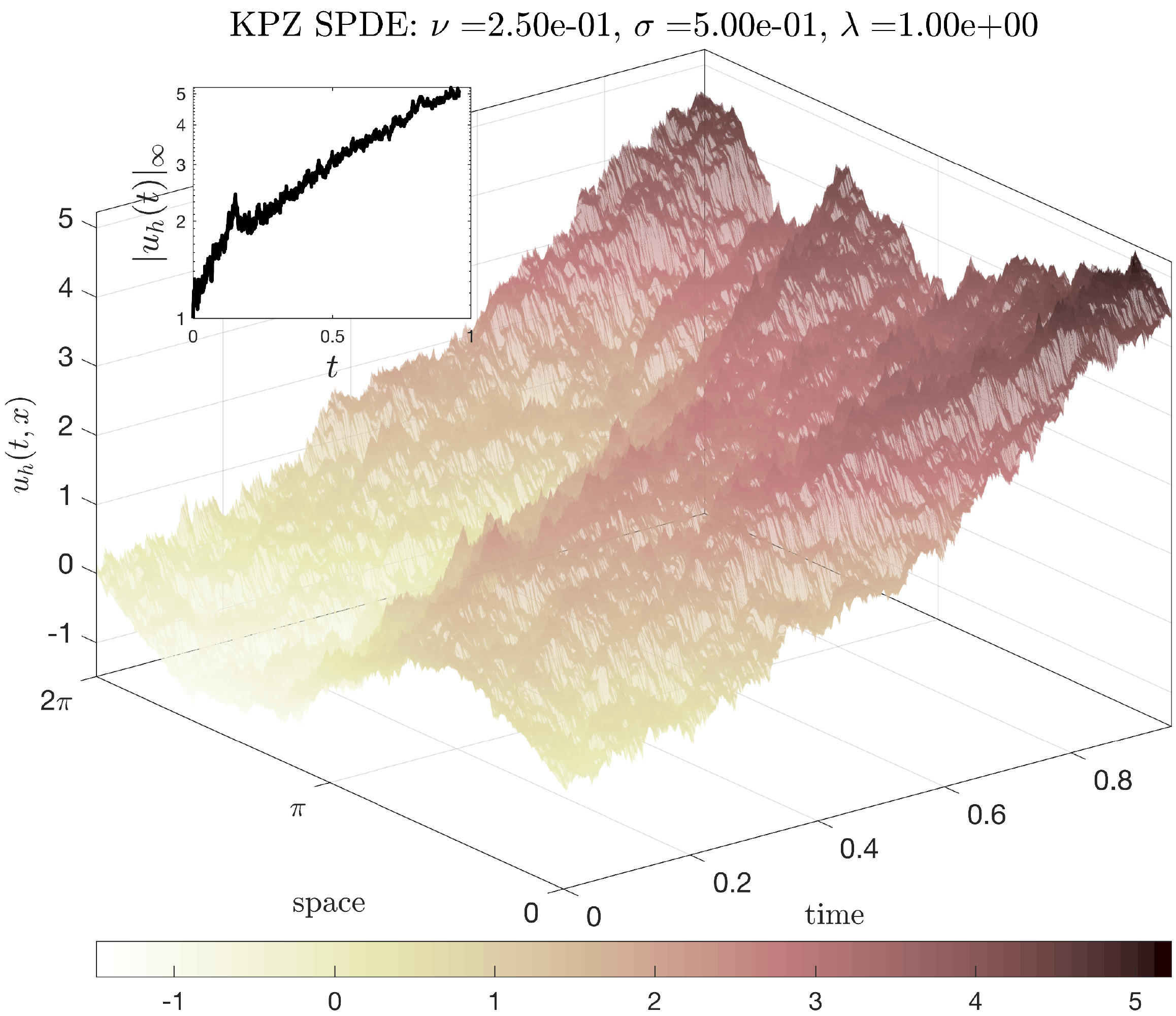}  \\
\includegraphics[width=0.65\textwidth]{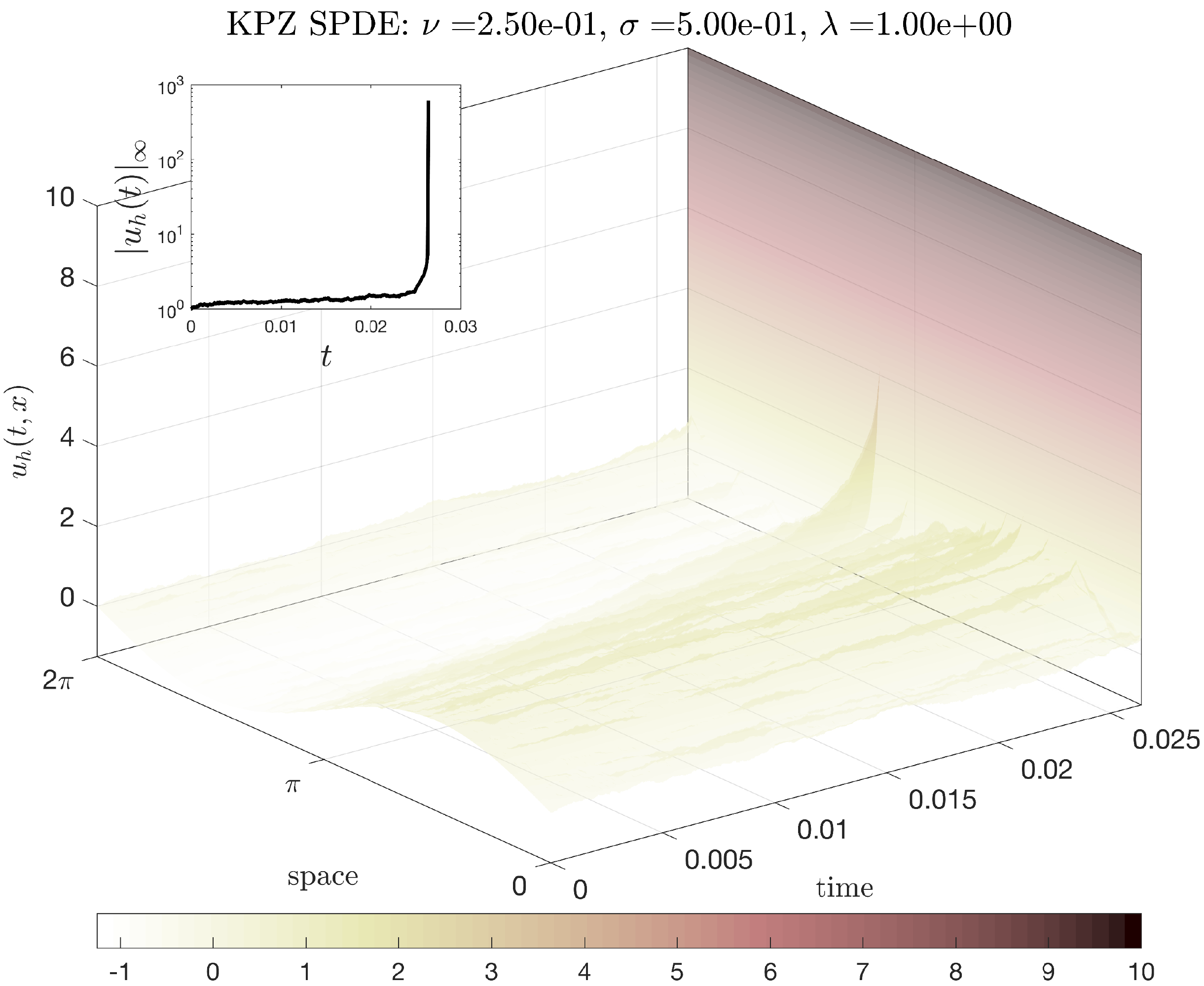} 
\end{center}
\caption{ \small {\bf KPZ SPDE on $[0,2 \pi]$.}   Sample paths produced using SPECTRWM at the parameter values indicated in the figure title.  Colors indicate the height of the surface, and are added in order to make the surface plots a bit more clear.  The left panel uses a valid semi-discrete approximation, while the right panel uses a semi-discrete approximation that converges to a non-ergodic SPDE.    }
 \label{fig:spectrwm_kpz_spde_sample_paths}
\end{figure}

Here we apply SPECTRWM to the following KPZ SPDE: \begin{equation}
\label{eq:heat_spde}
\begin{dcases}
du =  \left( \dfrac{\partial^2 u}{\partial x^2}  + \lambda \left( \frac{\partial u}{\partial x} \right)^2 \right) dt + \sigma d W & \forall x \in [0,2 \pi], t \ge 0 \\
u(t,0) =  u(t,2 \pi)  & \forall t\ge 0 \;, 
\end{dcases}
\end{equation}
with a sinusoidal initial condition as depicted in Figure~\ref{fig:spectrwm_kpz_spde_sample_paths}.
To construct approximating SDEs, we use the previous setup for the linear part and discretize the nonlinear term using either:\[
 \dfrac{\partial f}{\partial x}(x_i)^2  \approx \begin{dcases}
 \left( \frac{1}{2 \Delta x} \left( f_{(i+1) \bmod n}- f_{(i-1) \bmod n} \right) \right)^2 \\
 \left( \frac{1}{\Delta x} \left( f_{(i+1) \bmod n} - f_{i} \right)  \right)^2 
\end{dcases}
\]  

\section{Conclusion}

What this paper has achieved is a generalization of the 1D random walk approximation of Brownian motion to SPDEs.   Indeed, somewhat like the 1D random walk which makes jumps of fixed size $h$ to the left or right with equal probability, SPECTRWM jumps of fixed size $h$ forward or backward along the leading $n$ eigenfunctions of the linear part of the drift  with jump rates that depend only on its current state.  Moreover, just as the 1D random walk uses a mean holding time of $h^2$ for the sake of local accuracy, the mean holding time of SPECTRWM is also determined by local accuracy.  

Aside from its mathematical interest, this generalization has practical uses.  Indeed, since SPECTRWM is a (continuous-time) jump process with fixed jump size, this method solves in a natural way several problems that one encounters when simulating SPDEs, including their long-time simulation, the approximation of their path-dependent expected values, and by construction of its jumps, SPECTRWM is faithful to the domain of the SPDE solution.   Being a Markov jump process, SPECTRWM also has the advantage of being easy to quantitatively analyze since it admits an infinitesimal generator, which can be used to analyze its stability and accuracy \cite{BoVa2016, BoSa2015}.  

\appendix

\section{Local Consistency} \label{sec:local_consistency}

Here we show that the infinitesimal generator of SPECTRWM is an accurate approximation of the infinitesimal generator  of the approximating SDE given in \eqref{eq:semi_discrete}.  Referring to \eqref{eq:generator}, Taylor expand the terms $f(v \pm h e_i)$ to obtain: \begin{align*}
 Q f(v) = \sum_{1 \le i \le n} & ( J_i^+ - J_i^- ) \left( h \nabla f(v)^T e_i + \frac{h^3}{6} D^3 f(v)( e_i , e_i) + O(h^5) \right)  \\
  + & ( J_i^+ + J_i^- ) \left(  \frac{h^2}{2} D^2 f(v)( e_i , e_i) + \frac{h^4}{12} D^4 f(v)( e_i , e_i) + O(h^6)  \right)   
\end{align*}
To leading order, the forward/backward jump rates satisfy: \begin{align*}
( J_i^+ + J_i^- ) &= \frac{\sigma^2}{h^2 \Delta x} + O(\Delta x) \\ 
( J_i^+ - J_i^- ) &=  \frac{1}{h} ( \mu_i v^T e_i +  F_n(v)^T e_i ) + O( h \Delta x^2  ) 
\end{align*}
Thus, \begin{align*}
 Q f(v) = \sum_{1 \le i \le n} & ( \mu_i v^T e_i +  F_n(v)^T e_i ) \nabla f(v)^T e_i + \frac{\sigma^2}{2 \Delta x} D^2 f(v)( e_i , e_i) + O(h^2 \Delta x + h^2) 
\end{align*}
To finish, we use the fact that $e_i$ is an eigenvector of $L_n$ to obtain:  \begin{align*}
 Q f(v) &= \tr\left( \left( \nabla f(v)  v^T L_n +  \nabla f(v) F_n(v)^T + \frac{\sigma^2}{2 \Delta x} D^2 f(v) \right) \sum_{1 \le i \le n} e_i e_i^T \right) + O(h^2 \Delta x + h^2)   \\
 &=  v^T L_n \nabla f(v) + F_n(v)^T \nabla f(v) + \frac{\sigma^2}{2 \Delta x} \tr( D^2 f(v) ) + O(h^2 \Delta x + h^2) 
\end{align*}

\section*{Acknowledgements}

I wish to acknowledge Gideon Simpson for encouraging me to pursue this project, and 
the participants of the 2016 Gene Golub Summer School at Drexel University for their feedback.
I also wish to thank Eric Vanden-Eijnden for his helpful comments on an earlier version of this article.

\bibliography{nawaf}

\end{document}